\documentclass[aop,MSNbibl,dvips]{arximspdf}

\usepackage{graphicx}

\doi{10.1214/12-AOP747} 
\volume{41}
\issue{4}
\pubyear{2013}
\firstpage{2961}
\lastpage{2989}

\makeatletter
\newcommand{\eqref}[1]{(\ref{#1})}
\newcommand{\eps}{\varepsilon}
\newcommand{\E}{\mathsf{E}}
\newtheorem{lemma}{Lemma}[section]
\newcommand{\Fc}{\mathcal{F}}
\newtheorem{theorem}{Theorem}[section]

\newproclaim{remark}{Remark}[section]
\newcommand{\Pp}{\mathsf{P}}
\newcommand{\Ac}{\mathcal{A}}
\newcommand{\R}{\mathbb{R}}
\newcommand{\N}{\mathbb{N}}
\newcommand{\Z}{\mathbb{Z}}
\newcommand{\Sb}{\mathbb{S}}
\newcommand{\Nc}{\mathcal{N}}
\newcommand{\ONE}{\mathbf{1}}
\newcommand{\bg}{\bar{\gamma}}
\newcommand{\Ub}{\mathbb{U}}
\makeatother

\begin{document}
\begin{frontmatter}

\title{The Burgers equation with Poisson random~forcing}
\runtitle{Burgers equation with Poisson random forcing}

\begin{aug}
\author{\fnms{Yuri} \snm{Bakhtin}\corref{}\ead[label=e1]{bakhtin@math.gatech.edu}\thanksref{t1}}
\thankstext{t1}{Supported in part by NSF CAREER Grant DMS-07-42424.}
\runauthor{Y. Bakhtin}
\affiliation{Georgia Institute of Technology}
\address{School of Mathematics\\
Georgia Institute of Technology\\
686 Cherry Street\\
Atlanta, Georgia 30332-0160\\
USA\\
\printead{e1}}
\end{aug}

\received{\smonth{6} \syear{2011}}
\revised{\smonth{2} \syear{2012}}

%
\begin{abstract}
We consider the Burgers equation on the real line with forcing given by
Poissonian noise with no periodicity assumption.
Under
a weak concentration condition on the driving random force, we prove existence
and uniqueness of a global solution in a certain class. We describe its
basin of attraction that can also be viewed as
the main ergodic component for the model. We establish existence and
uniqueness of global minimizers associated to the
variational
principle underlying the dynamics. We also prove the diffusive behavior
of the global minimizers on the universal cover
in the periodic forcing case.
\end{abstract}

%


\end{frontmatter}

\section{Introduction}
The Burgers equation is one of the basic nonlinear evolution equations.
%
\begin{equation}
\label{eqburgers} \partial_t u(t,x)+u(t,x)\cdot\partial_x
u(t,x)= f(t,x).
\end{equation}
Here $t\in\R$ is the time variable, and $x\in\R$ is the space
variable. The equation describes the evolution of velocity
vector field $u(\cdot,\cdot)$ of sticky dust particles
in the presence of external potential forcing $f(t,x)=-\partial_x F(t,x)$.

Burgers introduced this equation as a turbulence model. Although it was
soon discovered that the dynamics governed
by~\eqref{eqburgers} does not describe turbulence adequately,
the equation has naturally appeared in various other contexts, from
cosmology to traffic modeling. An informative recent
survey
on Burgers turbulence is~\cite{Bec-KhaninMR2318582}.

One of the remarkable properties of the Burgers equation is that even
if the initial data at time~$t_0$
and the forcing are smooth, the solution of the Cauchy problem
typically develops discontinuities or shocks, and if one
wants to extend the solution
beyond the formation of shock waves, one has to work with generalized
solutions. Under mild assumptions on the initial
data and forcing, only one of
the generalized solutions is physical. This solution is called the
entropy or viscosity solution, and it can be found
using a
characterization that
is often called the Lax--Oleinik variational principle (see, e.g.,
\cite{Bec-KhaninMR2318582} and references therein).
Namely, the solution potential [a function $U$ such that $\partial_x
U(t,x)=u(t,x)$
for a.e. $x\in\R$] satisfies
%
\begin{equation}
U(t,x)=\inf_{\gamma: \gamma(t)=x} \biggl\{U\bigl(t_0,\gamma(t_0)
\bigr)+\int_{t_0}^t L\bigl(s,\gamma(s),\dot
\gamma(s)\bigr)\,ds \biggr\}. \label{eqvariationalprinciple}
\end{equation}
The expression in the curly brackets is called action, and the the
infimum of action
is taken over all absolutely continuous trajectories $\gamma$ defined
on $[t_0,t]$ and terminating at $x$ at time $t$.
The Lagrangian $L$ is defined by
\[
L(t,x,p)= \frac{p^2}{2}-F(t,x).
\]

Following the hydrodynamic interpretation of the Burgers equation, one
can identify the action minimizers
in~\eqref{eqvariationalprinciple} as the particle trajectories. This
kind of representation holds true for a more
general equation of Hamilton--Jacobi type.
The specifics of the Burgers equation is that if $\gamma^*$ is
a unique minimizer in \eqref{eqvariationalprinciple}, then
$u(t,x)=\dot\gamma^*(t)$.

When the forcing is a random field, one has to work with optimization
problems for paths accumulating
action from a random Lagrangian landscape, so questions about Burgers
equations with randomness
become random media questions.

The ergodic theory of the Burgers equation with random forcing begins
with~\cite{ekmsMR1779561}.
The forcing in~\cite{ekmsMR1779561} is assumed to be white noise type
in time and smooth and periodic in space. Due to
the periodicity assumption, the
evolution effectively takes place on a circle. The compactness of the
circle allows for efficient control of the long
time behavior of action minimizers, which leads
to constructing attracting global solutions and thus
to a complete description of the ergodic components for the dynamics,
each one consisting of all velocity profiles with
given mean velocity.

This work was extended and streamlined
in~\cite{IturriagaMR1952472} and \cite
{Gomes-Iturriaga-KhaninMR2241814}, where the multidimensional version
of the
Burgers equation with positive or zero
viscosity on a torus was considered. In~\cite{ybMR2299503} the
ergodic theory for the Burgers equation on a segment
with
random boundary conditions was developed.

In all of these papers the compactness of the domain played an
important role. In fact, in the case of unbounded domain
with no periodicity assumption, currently there is no complete
understanding of the ergodic properties of the Burgers
equation. Let us summarize what is known.

In~\cite{Khanin-HoangMR1975784}, the Burgers equation in $\R^d$ with
aperiodic white-noise forcing with certain
localization properties was considered.
A global solution constructed in the paper was shown to have a basin of
attraction containing the zero velocity profile,
but no interesting
properties of the global solution were established, and the description
of the domain of attraction of the global
solution was incomplete.

In~\cite{SuidanMR2141893}, it was noted that in the absence of
periodicity assumptions, the long time behavior of
solutions can
depend on the behavior of the initial condition at infinity in an
essential way. In particular, it was shown that
outside the main ergodic component (containing the
zero velocity profile), there are solutions with significantly
different behavior.

In this paper we introduce a new kind of random forcing for the Burgers
equation on the real line with no periodicity
assumption.
The forcing potential we suggest is given by a Poisson point field.
In this model, paths accumulate their action traveling through a
cloud of random Poissonian points. Although this model preserves many
features of the white noise model, it is easier to
analyze and visualize. It also has
much in common with the well-known Hammersley process (see, e.g.,~\cite
{Aldous-DiaconisMR1355056}) which has been
explicitly used for the analysis of hydrodynamic limit
resulting in the Burgers equation in~\cite{SeppMR1386297}.\looseness=1

For the new
model, we are able to construct a global solution via a limiting
procedure seeded at zero initial condition,
prove a so called one force---one solution principle (1F1S), and
describe the main ergodic component of the system,
that is, the
basin of attraction of the global solution.

1F1S for the Burgers equation on the circle is tightly connected to the
hyperbolicity of the global action minimizer.
In particular, for any two Burgers particles, the distance between
their backward trajectories (given by the
corresponding one-sided action minimizers)
converges to zero. A stronger phenomenon occurs in the case of
Poissonian forcing: for any two particles, their backward
trajectories will
meet at one of the Poissonian points in finite time and coincide from
that point on in the reverse time. This stronger
form of hyperbolicity
may naturally be called hyperhyperbolicity.

The rest of the paper is organized as follows:
In Section~\ref{secpoissonianmodel} we introduce the new forcing
model based on Poissonian points.
In Section~\ref{secgeometry-of-solutions} we discuss the geometry of
foliation of the space--time into particle
trajectories
under point forcing. In Section~\ref{secmainres} we formulate our
main results. In Section~\ref{secconstructing} we
construct
the global solution. In Section~\ref{secatinfty} we describe its
behavior at infinity. In Section~\ref{secattractor},
we
show that this solution is an attractor and describe its basin of
attraction. In Section~\ref{secglobalmin} we
study global minimizers. An important part of that section is a Central
limit theorem discribing the diffusive behavior
of global minimizers
for periodic Poissonian forcing.

\section{Poissonian point forcing}
\label{secpoissonianmodel}

The goal of this section is to describe the model rigorously, so let us
now be more precise. The model is based
on a Poisson point field, and we refer to~\cite{KallenbergMR854102}
for an introduction to point processes as random
integer-valued measures.

We are working on a complete probability space $(\Omega,\Fc,\Pp)$.
It is convenient to identify $\Omega$ with
the space of locally finite point configurations $\omega=\{(s_i,x_i),
i\in\N\}$ in space--time $\R\times\R$. The
sigma-algebra $\Fc$ is
generated by maps $N(B)$ assigning to each $\omega$ the number of
points of $\omega$ in a bounded Borel set
$B\subset\R\times\R$.
The measure~$\Pp$ is the distribution of a Poisson point field with
intensity measure $\mu(dt\times dx)$.

Since we want the forcing to be stationary in time, we shall always assume
that the intensity is a product measure
\[
\mu(dt\times dx)=dt\times m(dx).
\]
Then for disjoint sets $B_1,\ldots, B_n$, the random variables
$N(B_1),\ldots,N(B_n)$ are independent and Poissonian with parameters
$\mu(B_1),\ldots,\mu(B_n)$.

We will denote the integral term in~\eqref{eqvariationalprinciple} as
\[
\Ac^{t_0,t}(\gamma)=\int_{t_0}^t L\bigl(s,
\gamma(s),\dot\gamma(s)\bigr)\,ds=\frac{1}{2}\int_{t_0}^{t}
\dot\gamma^{2}(s)\,ds - \int_{t_0}^{t} F
\bigl(s,\gamma(s)\bigr)\,ds,
\]
and redefine the contribution from the potential by
%
\begin{equation}
\label{eqpoisson-forcing} \int_{t_0}^{t} F\bigl(s,
\gamma(s)\bigr)\,ds=N^{t_0,t}(\gamma),
\end{equation}
where for a path $\gamma$ and times $t_0$ and $t$ satisfying $t_0<t$,
$N^{t_0,t}(\gamma)=N^{t_0,t}_\omega(\gamma)$ is
the number of Poissonian points that $\gamma$ passes
through between $t_0$ and $t$. In other words, each Poissonian point
visited by the path contributes $-1$ to the
action.
An immediate generalization
of our model is a compound Poisson point field where each point comes
with a random weight which results in random
contributions to the action. In fact,
all our results can be extended to that case under reasonable
assumptions on the random weights. However, for simplicity
we concentrate here on the simple Poisson process.

Definition~\eqref{eqpoisson-forcing} results in the following
expression for
action accumulated by a path $\gamma$ between times $t_0$ and $t>t_0$:
\[
\Ac^{t_0,t_1}(\gamma)=\Ac_\omega^{t_0,t}(\gamma)=
\frac{1}{2}\int_{t_0}^{t}\dot
\gamma^{2}(s)\,ds - N^{t_0,t}_\omega(\gamma).
\]

It is well known (or can be easily derived from the Euler--Lagrange
equations) that in the zero forcing field, the
minimizers (or particle trajectories) are straight lines. We conclude
that between visits to Poissonian points, action-minimizing paths are
straight lines.

Let us introduce more notation. For two times $t_0$ and $t_1$, and two
sets $A_0,A_1\subset\R$, we denote by\vspace*{1pt}
$\Gamma_{t_0,A_0}^{t_1,A_1}$ the set of all piecewise linear paths
defined between $t_0$ and $t_1$ such that switchings
from one linear regime to another happen only at Poissonian points. We
also denote\vspace*{-1pt} the set of action minimizers over
$\Gamma_{t_0,A_0}^{t_1,A_1}$ by
$M_{t_0,A_0}^{t_1,A_1}=M_{t_0,A_0}^{t_1,A_1}(\omega)$.\vspace*{1pt} If $A_0$ or
$A_1$ consists of one
point~$x$, we will often use
index $x$ instead of $\{x\}$ in these notation.

To define the main random dynamical system, we must start with the
phase space. First we recall that the natural space of
solutions for the Burgers
equation consists of piecewise continuous functions $u$ defined on $\R
$, with right and left limits at every point, with
at most
countably many discontinuities, each discontinuity
being a downward jump or shock $u(x-)>u(x+)$. (The shock absorbs
incoming particles on both sides.)
We shall impose an additional restriction on these functions to be
bounded and will not distinguish between
two functions that coincide at all their continuity points. We will
denote the resulting factor space by~$\Ub$, and
often we will abuse the notation
writing $u\in\Ub$ when $u$ is a representative of an element of $\Ub$.

We will need a measure of proximity in $\Ub$. We denote the set of
continuity points of a function $h\in\Ub$ by $C_h$,
and for any
$h_1,h_2\in\Ub$, write
\[
d(h_1,h_2)=\exp\bigl[-\sup\bigl\{r>0\dvtx
h_1(x)=h_2(x), x\in C_{h_1}\cap C_{h_2}
\cap B_r \bigr\} \bigr],
\]
where $B_r=[-r,r]$. If there is no neighborhood of the origin where
$h_1$ and~$h_2$ coincide, we set $d(h_1,h_2)=1$.
If $h_1\equiv h_2$, we set $d(h_1,h_2)=0$.
Thus defined~$d$ is a metric in $\Ub$ taking values in $[0,1]$.

Given $v\in\Ub$, we can define a potential $V$ so that $V'(x)=v(x)$
for all $x$. For any times $t_0,t_1$ with $t_0<t_1$,
we set
%
\begin{equation}
\label{eqdefofcocycle} \Phi^{t_0,t_1}_\omega v(x)=\dot
\gamma^*(t_1),
\end{equation}
where $\gamma^*$ is the solution of
%
\begin{equation}
\label{eqminimizationproblem} A_V^{t_0,t_1}(\gamma)=V\bigl(
\gamma(t_0)\bigr)+\Ac^{t_0,t_1}_\omega(\gamma) \to\min,\qquad
\gamma\in\Gamma_{t_0,\R}^{t_1,x}.
\end{equation}
Let us assume that
%
\begin{equation}
\label{eqfiniteness-of-poisson} m(\R)<\infty,
\end{equation}
and briefly summarize (without proof) several facts about the Burgers
equation solution map $\Phi$ that apply to the
current setting.

\begin{lemma} If $h\in\Ub$, then with probability 1 the following hold:
\begin{longlist}[(1)]
\item[(1)] For any time interval $[t_0,t_1]$, in definition \eqref
{eqdefofcocycle}--\eqref{eqminimizationproblem},
the minimizer~$\gamma^*$ [and, consequently, its slope $\dot\gamma
^*(t_1)$ at the terminal time $t_1$] is
defined uniquely for all $x\in\R$ except at most countably many points.
Every point $x$ where $\Phi^{t_0,t_1}_\omega v(x)$ is uniquely defined
is a continuity point of $\Phi^{t_0,t_1}_\omega v$. At any point where
the minimizer is not unique,
$\Phi^{t_0,t_1}_\omega v$ makes a downward jump.
\item[(2)] The function $\Phi^{t_0,t_1}_\omega v$ is bounded (in
particular, combining this with the first part of this
lemma, we obtain
that $\Phi^{t_0,t_1}_\omega$ is a map from $\Ub$ to itself).
\item[(3)] Moreover, for all $\omega$, if $t_0\le t_1\le t_2$,
%
\begin{equation}
\label{eqcocycle1} \Phi^{t_1,t_2}_\omega\Phi^{t_0,t_1}_\omega
v= \Phi^{t_0,t_2}_\omega v.
\end{equation}
\end{longlist}
\end{lemma}
%
\begin{remark}
Introducing $\Phi^{t}_\omega=\Phi^{0,t}_\omega$ for $t\ge0$, we
can rewrite the cocycle property~\eqref{eqcocycle1} as
\[
\Phi^{t_1+t_2}_\omega v=\Phi^{t_2}_{\theta^{t_1}\omega}
\Phi^{t_1}_\omega v,\qquad t_1,t_2\ge0,
\]
where $\theta^t$ denotes the time shift of the Poissonian point field
$(s_i,x_i)\mapsto\break (s_i-t,x_i)$.
\end{remark}

Let us denote by $\Fc_A$ the sigma-algebra generated by the
restriction of the Poissonian point field to $A\times\R$ for
any set $A$ of times.
Clearly, the random operator $\Phi_\omega^{t_0,t_1}$ depends only on
the realization of the Poisson process between
times $t_0$ and~$t_1$; that is, it is measurable w.r.t. $\Fc_{[t_0,t_1]}$.

\section{Geometry of solutions under Poissonian forcing}\label
{secgeometry-of-solutions}

Throughout this paper, we consider the external forcing that is
concentrated on a discrete set of Poissonian points (we
will often call
them \textit{forcing points}). This is different from
traditionally considered smooth forcing fields, so let us understand
the effect of this kind of forcing on the
solution.

Let us consider a model situation where a smooth beam of Burgers
particles encounters a forcing point at the origin at
time $0$.
Let us assume that at time~$0$, the velocity vector field near $0$ is
$u_0(y)=a+by$, where $b>0$.

It is clear that
for every $(t,x)$ with $t>0$ and $x$ close to the origin, there are two
minimizer candidates. The minimizer either
passes through the origin, or it does not.
If it does, then (assuming there are no other point sources of forcing)
it has to be a straight line connecting the
origin to $(t,x)$,
and the accumulated action is $A_1(t,x)=x^2/(2t)-1$, where $-1$ is the
contribution of the forcing point at the origin,
and $x^2/(2t)$
is the action accumulated while moving with constant velocity $x/t$
between $0$ and $t$. If the minimizer
does not pass through the origin, then it is a straight line connecting
some point $(0,x_0)$ to $(t,x)$. On the one
hand,
the velocity of the particle associated with the minimizer is
$(x-x_0)/t$. On the other hand, it has to coincide
with $u_0(x_0)=a+bx_0$. Therefore, we can find $x_0=(x-at)/(1+bt)$.
Taking into account that $U_0(x_0)=ax_0+bx_0^2/2$,
we can compute that the total action of that path is
$A_2(t,x)=(bx^2+2ax-a^2t)/(2(1+bt))$.

\begin{figure}

\includegraphics{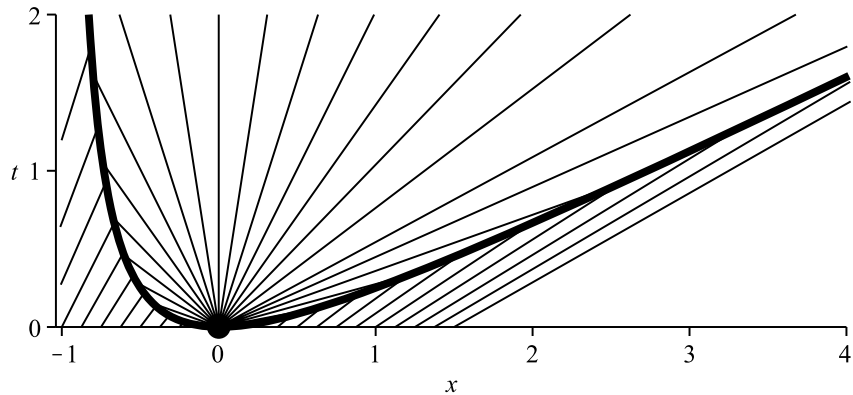}

\caption{Minimizers around a forcing point.}\label{fig1pt}
\end{figure}

To see which of the two cases is realized for $(t,x)$ we must compare
$A_1(t,x)$ and $A_2(t,x)$. If $A_1(t,x)<A_2(t,x)$,
then the particle
arriving to $x$ at time $t$ is at the origin at time $0$. If
$A_1(t,x)>A_2(t,x)$, then the particle arriving to $x$ at
time~$t$ is
one of the particles that moved with constant velocity and was a part
of the incoming beam. If $A_1(t,x)=A_2(t,x)$,
then
both of these paths are minimizers, and at time $t$ there is a shock at
point $x$. The relation $A_1(t,x)=A_2(t,x)$ can
be rewritten
as
\[
(x-at)^2=2t(1+bt).
\]
For small values of $t$, the set of points satisfying this relation
looks like a parabola $(x-at)^2=2t$; see
Figure~\ref{fig1pt} where an example with $a=1$ and $b=1/2$
is shown.\vadjust{\goodbreak}

We see that when a Poissonian point appears, it emits a continuum of
particles, each moving with constant velocity,
creating two shock fronts moving (at least for a short time) to the
left and right.

It is important to notice that in our model case with a forcing point
at the origin, $u(t,x)=x/t$ for all points
connected to the origin by a minimizing segment. It means that
for each time $t$, the velocity is linear in the domain of influence of
the forcing point, and the velocity gradient
decays with time
as $1/t$.

In general, the behavior of this kind occurs near each forcing point,
and in the long run more and more points of the
space--time plane get
assigned to forcing points. Grouping together points assigned to the
same forcing point, we obtain a tesselation of
space--time
into domains of influence of forcing points. Inside each domain or
cell, the velocity field is linear in $x$ if the
time~$t$ is fixed.

It is well known that in the Burgers equation the energy is dissipated
at the shocks; see, e.g.,
\cite{Bec-KhaninMR2318582}.
By seeding new particles at each Poissonian point,
the forcing pumps energy into the system, and, therefore, we can hope
that there is a dynamical or statistical energy
balance in the system.
We will actually see that this dissipation results in asymptotic
alignment of the velocities of particles that keep
moving away from the origin without
being absorbed into shocks.

Another point of view at the stationarity and ergodicity issues for
this system is related to the stabilization of the
tesselation
of space--time into cells described above.

\section{Main result}\label{secmainres}
Although
it would be interesting to consider the situation where the spatial
intensity measure $m(dx)$ satisfies $m(\R)=\infty$
(e.g., the Lebesgue measure on $\R$),
throughout this paper we will adopt either assumption~\eqref
{eqfiniteness-of-poisson} or an even stronger finite first
moment assumption
%
\begin{equation}
\int_\R\bigl(1+|x|\bigr)m(dx)< \infty. \label{eqmomentcondition}
\end{equation}

\begin{theorem} \label{thmain} Suppose~\eqref{eqmomentcondition}
holds. Then there is a set $\Omega'$ with
$\Pp(\Omega')=1$ and
a function $u\dvtx\R\times\Omega'\to\Ub$ such that on $\Omega'$ the
following hold:
\begin{longlist}
\item[(1)] \hypertarget{itmeasurability-of-u} $u$ is measurable w.r.t.
$\Fc_{(-\infty,0]}$. In other words, it depends only on
$\omega|_{(-\infty,0]}$.
\item[(2)]\hypertarget{itglobal-solution} $u$ defines a global solution
[in other words, it is skew-invariant under
$(\Phi,\theta)$]
\[
\Phi^t_\omega u_\omega=u_{\theta^t\omega},\qquad t\ge0.
\]
\item[(3)] \hypertarget{itpiecewise-linear}The solution $u_\omega$ is
piecewise linear.
\item[(4)] \hypertarget{itrho} There is a nonrandom constant $q>0$ such that
\[
\lim_{x\to\pm\infty}u_\omega(x)=\pm q.
\]
\item[(5)] \hypertarget{itattractor}
This solution $u$ plays the role of a one-point attractor. Namely, if
$V'=v\in\Ub$ and
%
\begin{equation}
\label{eqasymptoticslope} \liminf_{x\to\infty} \frac{V(x)}{x}> -q,
\end{equation}
we have forward attraction,
%
\begin{equation}
\label{eqfwdattr} d\bigl(\Phi^t_\omega v, u_{\theta^t\omega}
\bigr)\to0,\qquad t\to\infty,
\end{equation}
and pullback attraction,
%
\begin{equation}
\label{eqpullattr} d\bigl(\Phi^t_{\theta^{-t}\omega} v,  u_\omega
\bigr)\to0,\qquad t\to\infty.
\end{equation}
\end{longlist}
The function $u$ is a unique (up to zero measure modifications) global
solution satisfying
%
\begin{equation}
\label{equniqunessclass} \liminf_{x\to\infty} \frac{U_\omega
(x)}{x}>- q,
\end{equation}
with positive probability (here $U_\omega$ is the potential of
$u_\omega$, i.e., $U'_\omega\equiv u_\omega$).
\end{theorem}

\begin{remark} One can reformulate the theorem in terms of a global
solution defined as a function of three
variables,
$u_\omega(t,x)=u_{\theta^t\omega}(x)$.
\end{remark}

\begin{remark}
If one accepts a weaker condition~\eqref{eqfiniteness-of-poisson},
then all conclusions
of Theorem~\ref{thmain} except conclusion~(\hyperlink{itrho}{4}) still hold,
and their proofs do not change.
Conclusion~(\hyperlink{itrho}{4}) has to be replaced
with a weaker one,
\[
\lim_{x\to\infty} \frac{U_\omega(x)}{x}=q.
\]
\end{remark}

\begin{remark} Conclusion~(\hyperlink{itrho}{4}) means that in the stationary
regime, at infinity one observes particles moving
away from the origin with velocity $q$.\vadjust{\goodbreak}
\end{remark}

\begin{remark} Conclusion~(\hyperlink{itattractor}{5}) means that if one
starts with an initial condition that sends particles
from infinity toward zero
with speed that is less than $q$ [see condition~\eqref
{eqasymptoticslope}], then this inbound flow is not strong enough
to compete with the outbound flow of particles
developed due to the noise, and in the long run it is dominated by the
latter. If
condition~\eqref{eqasymptoticslope} is violated, then the long term
properties of solutions are sensitive to the details of the behavior of
the initial condition at infinity because the
inbound flow of particles
may be stronger than the outbound one, and one will observe effects
similar to those discussed
in~\cite{SuidanMR2141893}.
\end{remark}

\begin{remark} The uniqueness conclusion~(\hyperlink{itattractor}{5}) and
measurability property [conclusion
(\hyperlink{itmeasurability-of-u}{1})] can be
combined into 1F1S principle---at time 0 there is a unique velocity
profile compatible with the history of the
forcing.
\end{remark}

\section{Constructing a global solution}\label{secconstructing}

In this section we construct a global solution $u$. To do this, we
start with the zero initial condition at time $-T$
and take $T$ to infinity. Our goal is to show that $\Phi^{-T,0}_\omega
0$ converges in $(\Ub,d)$ to a limiting function
and that
this limit defines a global solution.

It will be convenient to assume that 0 belongs to the support of
measure~$m$, i.e., for any $\delta>0$,
$m(B_\delta)>0$.
We adopt this nonresrictive assumption without loss of generality since
one can always introduce a shift coordinate
change to make it hold true.

Since all admissible paths are composed of straight line segments, we
will often use the following elementary result on
action accumulated along one segment:

\begin{lemma}\label{lmaction-for-segment} A path corresponding to a
particle moving with constant velocity $v$ for
time $t$ and visiting no Poissonian points,
accumulates action equal to
\[
\frac{v^2t}{2}=\frac{vx}{2}=\frac{x^2}{2t},
\]
where $x=vt$ is the traveled distance.
\end{lemma}

\begin{lemma}
\label{lmquasi-optimalpath}
There are numbers $a,b>0$ and an a.s.-finite random variable $\beta>0$
such that if $t>0$, $x\in\R$, and
$\omega\in\Omega$ sastisfy $t-|x|-2b\ge\beta(\omega)$, then
there is a path $\bar\gamma$ with $\bar\gamma(-t)=x$ and
\[
\Ac^{-t,0}_\omega(\bar\gamma)<-\bigl(t-|x|\bigr)a+ |x|+b.
\]
\end{lemma}
\begin{pf}Let us consider sets $A_k=[-2k,-2k+1]\times
B_{1/2}$. For any $k\in\N
$, we have
\[
\Pp\bigl\{N(A_k)\ne0\bigr\}=1 - e^{-M},
\]
where $N(A_k)$ denotes the number of Poisson points in $A_k$ and $M=m(B_{1/2})$.
For any $s>0$ we denote by $X(s)$ the random number\vadjust{\goodbreak} of indices $k\in\N
$ satisfying $k<s/2$ and $N(A_k)\ne0$.
The sequence $\ONE_{\{N(A_k)\ne0\}}$ is i.i.d. with mean $1-e^{-M}$,
and the strong law of large numbers implies that there is a random time
$\beta$ such that if $s>\beta$, then
%
\begin{equation}
\label{eqllnforR} X(s)>s\bigl(1-e^{-M}\bigr)/3.
\end{equation}
Consider a path $\bar\gamma$ that starts at $(-t,x)$ and
visits exactly one point in set~$A_k$ if $k$ satisfies $N(A_k)\ne0$ and
$2k<t-|x|-1$, and no other points.
Each Poissonian point in the path contributes $-1$ to the action, and
we can use Lemma~\ref{lmaction-for-segment} to
see that
each segment connecting these points
contributes at most $(2\cdot1/2)^2/2=1/2$. The slope of the segment
with endpoint $(-t,x)$
does not exceed $1$, and contributes at most $(|x|+1)/2$ to the action.

If $t-|x|-1>\beta$, then we can combine this with~\eqref
{eqllnforR} applied to $t-|x|-1$ and obtain
\begin{eqnarray*}
\Ac_\omega^{-t,0}(\bar\gamma)&<& - \bigl(t-|x|-1\bigr)\frac{1-e^{-M}}{3}
\frac
{1}{2}+\bigl(|x|+1\bigr)/2
\\
&<& - \bigl(t-|x|\bigr) \bigl(1-e^{-M}\bigr)/6 +|x|+\bigl(1-e^{-M}
\bigr)/6+1/2,
\end{eqnarray*}
and the lemma follows with $a=(1-e^{-M})/6$ and
$b=(1-e^{-M})/6+1/2>1/2$ since $t-|x|-2b>\beta$ implies
$t-|x|-1>\beta$.
\end{pf}

Let us recall that $B_r=[-r,r]$ for any $r>0$.
The following is the main localization lemma.
%
\begin{lemma}\label{lmmain-loc}
There are random variables $r^-$, $r^+$, $r^{\pm}$, $(\tau
^-_R)_{R>0}$, $(\tau^+_R)_{R>0}$, $(\tau^\pm_R)_{R>0}$,
such that for any $R>0$:
%
\begin{eqnarray}
\label{eqr-}
\qquad &&\Pp\bigl\{\mbox{there are } t>\tau^-_R, x\in
B_R \mbox{ and } \gamma\in M_{-t,x}^{0,\R} \mbox{ s.t. } \bigl|\gamma(0)\bigr|>r^- \bigr\}=0;
\\
\label{eqr+}
&&\Pp\bigl\{\mbox{there are } t>\tau^+_R, x\in
B_R \mbox{ and } \gamma\in M^{t,x}_{0,\R} \mbox{ s.t. } \bigl|\gamma(0)\bigr|>r^+ \bigr\}=0;
\\
\label{eqtwo-sided-r1}
&&\Pp\bigl\{\mbox{there are } t_-,t_+>\tau
^{\pm}_R,
x_-,x_+\in B_R \mbox{ and } \gamma\in M_{-t_-,x_-}^{t_+,x_+}
\mbox{ s.t. } \bigl|\gamma(0)\bigr|>r^\pm\bigr\}
\nonumber
\\[-8pt]
\\[-8pt]
\nonumber
&&\qquad=0.
\end{eqnarray}

Additionally, there are random variables $(\bar\tau_R)_{R>0}$ and a
number $R'>0$ such that for any $R>R'$,
%
\begin{eqnarray}
\label{eqtwo-sided-r2}\qquad &&\Pp\bigl\{\mbox{there are } t_->\bar\tau_R ,
t_+>\tau^{\pm
}_R, x_+\in B_R \mbox{ and }
\gamma\in M_{-t_-,\R}^{t_+,x_+} \mbox{ s.t. } \bigl|\gamma(0)\bigr|>r^\pm
\bigr\}
\nonumber
\\[-8pt]
\\[-8pt]
\nonumber
&&\qquad=0.
\end{eqnarray}
\end{lemma}
%
\begin{remark} The idea of this lemma is that minimizers over long
time intervals are localized within a random
neighborhood of the origin. Each of the random
variables $r^-,r^+$ and $r^{\pm}$ can be called localization radius.
\end{remark}

\begin{pf}
Let us prove \eqref{eqr-} first.
We are going to construct random variables~$K$ and $h$ so that for any
$R>0$ and for sufficiently large $t$,\vspace*{1pt} no path
$\gamma$ with $|\gamma(0)|>Kh$ can belong to $M_{-t,x}^{0,\R}$ with
$x\in B_R$. The reason why we need two random
variables is that we will use $h$ as an intermediate threshold.

Let $x\in B_R$. Consider a path $\gamma$ defined on $[-t,0]$ such that
$|\gamma(0)|>Kh$ and $\gamma(-t)=x$. Suppose
that $|\gamma(-s)|\le h$ for some $s\in[0,t]$ and define
$\sigma=\sup\{s\le t\dvtx |\gamma(-s)|\ge h\}.$ Then
\[
\Ac^{-\sigma,0}_\omega(\gamma)\ge\frac{(K-1)^2h^2}{2\sigma}-
N\bigl([-
\sigma,0]\times B_{h}^c\bigr).
\]
To treat the second term on the right-hand side, we need the following result:

\begin{lemma}
\label{lmlln-bound}
For any $\eps>0$, there is a positive random variable $R_0$ such that
with probability 1, for every $t>0$,
\[
N\bigl([-t,0]\times B_{R_0}^c\bigr) < \eps t.
\]
\end{lemma}

\begin{pf}Let us choose a number $\alpha_1$ such that
$m(B_{\alpha
_1}^c)<\eps/2$. Due to the strong law of large numbers,
there is a random time
$\tau>0$ such that $N([-t,0]\times B_{\alpha_1}^c)<\eps t$ for all
$t>\tau$. With probability 1, there are finitely many
Poissonian points
in $[-\tau,0]\times\R$. Let $\alpha_2$ be the maximal absolute
value of the spatial components of these points. The
conclusion of
the lemma holds true with $R_0=\alpha_1\vee\alpha_2$.
\end{pf}

Coming back to the proof of~\eqref{eqr-}, let us set $\eps=a/2$,
where $a$ is defined in
Lemma~\ref{lmquasi-optimalpath}.
Lemma~\ref{lmlln-bound} applied to this value of $\eps$
ensures the existence of $h=h(\omega)$ such that
%
\begin{equation}
\Ac^{-\sigma,0}_\omega(\gamma)\ge\frac{(K-1)^2h^2}{2\sigma} -
\eps\sigma.
\label{eqlowerbound2}
\end{equation}
If $\sigma<(K-1)h/\sqrt{2\eps}$, then $A_\omega^{-\sigma,0}(\gamma
)>0$, and the comparison with a zero velocity
trajectory
with zero action proves that $\gamma$ cannot be a minimizer.

To treat the case where
%
\begin{equation}
\label{eqthecaseoflargesigma1} \sigma\ge(K-1)h/\sqrt{2\eps},
\end{equation}
we will impose some restrictions on $K$. First, we require that $
K(\omega)\ge K_1(\omega)$, where
\[
K_1(\omega)= \biggl(\frac{\beta(\omega)+2b}{h(\omega)}+1 \biggr
)\sqrt{2\eps}+2,
\]
with $\beta$ and $b$ constructed in Lemma~\ref{lmquasi-optimalpath}.

Then, under assumption~\eqref{eqthecaseoflargesigma1}, $\sigma
-h-2b>\beta$ and we can apply
Lemma~\ref{lmquasi-optimalpath}. The path $\bar\gamma$ constructed
in that lemma for point $(-\sigma,\gamma(-\sigma))$
satisfies
\[
\Ac^{-\sigma,0}_\omega(\bar\gamma)\le-(\sigma-h)a+h+b.
\]
On the other hand, \eqref{eqlowerbound2} implies
\[
\Ac_\omega^{-\sigma,0}(\gamma)\ge- \eps\sigma\ge-a\sigma/2.
\]
The last two inequalities imply
\begin{eqnarray*}
\Ac^{-\sigma,0}_\omega(\bar\gamma)-\Ac^{-\sigma,0}_\omega
(\gamma)& \le&-(\sigma-h)a+h+b+ a\sigma/2
\\
&\le&-a\sigma/2 +h(a+1)+b,
\end{eqnarray*}
and, due to \eqref{eqthecaseoflargesigma1}, the right-hand side is
negative if we assume that $K(\omega)\ge K_2(\omega)$,
where
\[
K_2(\omega)=\frac{2\sqrt{2\eps}(h(\omega)(a+1)+b)}{ah(\omega)}+2.
\]
Therefore, under this assumption $\gamma$ cannot be a minimizer.
We conclude that if $K>K_1\vee K_2$, then $\gamma$ with
$|\gamma(0)|>Kh$ cannot be a minimizer satisfying $|\gamma(-s)|\le h$
for some $s\in[0,t]$.

Let us now consider a path $\gamma$ with $|\gamma(-s)|>h$ for all
$s\in[0,t]$.
We have then
\[
\Ac_\omega^{-t,0}(\gamma)\ge- \eps t\ge-at/2.
\]
On the other hand, we can invoke Lemma~\ref{lmquasi-optimalpath} to
see that
if $t-|x|-2b>\beta$, then
there is a path $\bar\gamma$ with $\bar\gamma(-t)=x$ such that
\[
\Ac_\omega^{-t,0}(\bar\gamma)<-\bigl(t-|x|\bigr)a+|x|+b
\]
and
\begin{eqnarray*}
\Ac_\omega^{-t,0}(\bar\gamma)-\Ac_\omega^{-t,0}(
\gamma)& \le&-\bigl(t-|x|\bigr)a+|x|+b+ at/2
\\
&\le&-at/2 +(a+1)|x|+b.
\end{eqnarray*}
Since $|x|\le R$, the right-hand side is negative if we require that
\[
t>\frac{2(R(a+1)+b)}{a}.
\]
Under this additional assumption, $\gamma$ cannot be a minimizer.
We conclude that \eqref{eqr-} holds if one chooses
\[
r^-(\omega)= \bigl(K_1(\omega)\vee K_2(\omega)\bigr)h(
\omega)
\]
and
\[
\tau^-_R(\omega)= \frac{2(R(a+1)+b)}{a} \vee\bigl(\beta(\omega)+R+2b
\bigr).\vspace*{1pt}
\]

The second part of the lemma, equation~\eqref{eqr+}, is only a time
reversed version of the first one.
The proof of~\eqref{eqtwo-sided-r1} is an adaptation of the above
argument to the two-sided situation.

Let us prove~\eqref{eqtwo-sided-r2}. First, choose $R'$ large enough
to ensure that due to the law
of large numbers, an optimal path cannot stay\vspace*{1pt} infinitely
outside~$B_{R'}$. Therefore, for sufficiently large $t_-$,\vspace*{-1pt}
minimizers from $M_{-t_-,\R}^{t_+,x_+}$ visit a point $x_-\in
B_{R'}\subset B_R$ between $-t_-$ and $-\tau_{R}^\pm$.
Since $x_-,x_+\in B_R$
and a restriction of a minimizer is a minimizer itself, we can finish
the proof by invoking~\eqref{eqtwo-sided-r1}.
\end{pf}

Let us denote
\begin{eqnarray*}
r(\omega)&=&r^-(\omega)\vee r^+(\omega)\vee r^{\pm}(\omega),
\\
\tau_R(\omega)&=&\tau_R^-(\omega)\vee
\tau_R^+(\omega)\vee\tau_R^{\pm}(\omega),\qquad R>0,
\\
D_1(R,T)&=&\bigl\{r(\omega)<R,r\bigl(\theta^T\omega
\bigr)<R, \tau_{R}(\omega)<T, \tau_{R}\bigl(
\theta^T\omega\bigr)<T\bigr\},\qquad R,T>0.
\end{eqnarray*}

\begin{lemma}\label{lmloc-with-positive-prob}
For any $L>0$ there are numbers $R>L$ and $T>0$ such that $P(D_1(R,T))>0$.
\end{lemma}

\begin{pf}We take $R$ so large that $\Pp\{r(\omega)>R\}
<1/4$. Then $\Pp\{
r(\theta^t\omega)>R\}<1/4$ for any $t$ since
$\theta^t$ preserves the measure. Then we take $T$ so large that
$\Pp\{\tau_{R}(\omega)>T\}<1/4$. Then $\Pp\{\tau_{R}(\theta
^T\omega)>T\}<1/4$, and the lemma follows.
\end{pf}

Let us fix the values of $R$ and $T$ given by Lemma~\ref
{lmloc-with-positive-prob} and introduce a new event
$D_2(R,T)$ consisting of all outcomes $\omega$ admitting
a point $(t^*,x^*)=(t^*,x^*)(\omega)\in[0,T]\times\R$ such that
for any $x,y\in B_R$, the
optimal path connecting $(0,x)$ and $(T,y)$ passes through $(t^*,x^*)$.

\begin{lemma}
\label{lmcollapse-with-positive-prob} Let $R$ and $T$ be provided by
Lemma~\ref{lmloc-with-positive-prob}. Then
\[
\Pp\bigl(D_1(R,T)\cap D_2(R,T)\bigr)>0.
\]
\end{lemma}

\begin{pf}The proof of this lemma is based on a
resampling of the point
configurations in $[0,T]\times\R$ according to a certain
kernel.
In this proof it is convenient to represent $\omega\in\Omega$ as
$\omega=(\omega_{\mathrm{in}},\omega_{\mathrm{out}})$ where
$\omega_{\mathrm{in}}\in\Omega_{\mathrm{in}}$ and $\omega_{\mathrm{out}}\in\Omega_{\mathrm{out}}$ are
restrictions
of the point configuration $\omega$ to $[0,T]\times\R$ and its
complement. We also denote by $\Pp_{\mathrm{in}}$ and $\Pp_{\mathrm{out}}$
the distributions of
Poisson point field in $[0,T]\times\R$ and its complement in $\R
\times\R$.

We will take a large number $n$ and consider a family of rectangles
$L_k,k=1,\ldots,n$ in $[0,T]\times\R$. We postpone a
precise
description of these rectangles.

For every $\omega\in D_1=D_1(R,T)$ we consider a new random
configuration $\omega'$. It coincides with
$\omega$ outside of $[0,T]\times\R$, and the restriction of $\omega
'$ onto $[0,T]\times\R$ consists of $n$
independent random points such that for each $k=1,\ldots,n$, the
distribution of $k$th point is concentrated in $L_k$,
$k=1,\ldots,n$.
Let us denote the distribution of the configuration\vadjust{\goodbreak} of these $n$ points
in $[0,T]\times\R$ by $\Pp'_{\mathrm{in}}$.
Later, we shall choose the distributions of individual points
appropriately to make $\Pp'_{\mathrm{in}}$ absolutely continuous
w.r.t. $\Pp_{\mathrm{in}}$.

To define the resampling more formally, for any $\omega$ we consider a
version of conditional probability
$\Pp(\cdot|\omega_{\mathrm{out}})$ defined for a set $D$ by
\[
\Pp(D|\omega_{\mathrm{out}})= \Pp_{\mathrm{in}}\bigl\{\omega'_{\mathrm{in}}
\dvtx\bigl(\omega'_{\mathrm{in}},\omega_{\mathrm{out}}\bigr)\in D
\bigr\},
\]
and define a new measure $\Pp'$ via
\[
\Pp'(E|\omega_{\mathrm{out}})=\Pp(D_1|
\omega_{\mathrm{out}})\Pp'_{\mathrm{in}}\bigl\{
\omega'_{\mathrm{in}}\dvtx\bigl(\omega'_{\mathrm{in}},
\omega_{\mathrm{out}}\bigr)\in E\bigr\}
\]
and
\[
\Pp'(E)=\int_{\Omega_{\mathrm{out}}} \Pp_{\mathrm{out}}(d
\omega_{\mathrm{out}})\Pp'(E|\omega_{\mathrm{out}}).
\]

Let us prove that $\Pp'\ll\Pp$. We must show that for any set $E$
with $\Pp'(E)>0$, we have $\Pp(E)>0$. Since
$\Pp'(E)>0$, the definition of $\Pp'$ yields
%
\begin{equation}
\Pp_{\mathrm{out}}\bigl\{\omega_{\mathrm{out}}\dvtx\Pp'(E|
\omega_{\mathrm{out}})>0\bigr\}>0. \label{eqpositiveset}
\end{equation}
Notice that if $\omega_{\mathrm{out}}$ satisfies $\Pp'(E|\omega_{\mathrm{out}})>0$,
then $\Pp(D_1|\omega_{\mathrm{out}})>0$ and
$\Pp'_{\mathrm{in}}\{\omega'_{\mathrm{in}}\dvtx\break  (\omega'_{\mathrm{in}},\omega_{\mathrm{out}})\in E\}
>0$. The latter
and the absolute continuity of $\Pp'_{\mathrm{in}}$ w.r.t. $\Pp_{\mathrm{in}}$ imply
that $\Pp_{\mathrm{in}}\{\omega'_{\mathrm{in}}\dvtx(\omega'_{\mathrm{in}},\omega_{\mathrm{out}})\in E\}>0$
for such $\omega_{\mathrm{out}}$. Therefore, due
to~\eqref{eqpositiveset},
\[
\Pp(E)= \Pp_{\mathrm{out}}\times\Pp_{\mathrm{in}}(E)=\int
_{\Omega} \Pp_{\mathrm{out}}(d\omega_{\mathrm{out}})
\Pp_{\mathrm{in}}\bigl\{\omega'_{\mathrm{in}}\dvtx\bigl(
\omega'_{\mathrm{in}},\omega_{\mathrm{out}}\bigr)\in E\bigr\} > 0,
\]
and the absolute continuity is proven.

Therefore, $\Pp(D_1(R,T)\cap D_2(R,T))>0$ will hold if
%
\begin{equation}
\label{eqp-prime-positive} \Pp'\bigl(D_1(R,T)\cap
D_2(R,T)\bigr)>0.
\end{equation}
So, it remains to finish the construction of the measure $\Pp'_{\mathrm{in}}$
and ensure that \eqref{eqp-prime-positive} holds
along with $\Pp'_{\mathrm{in}}\ll\Pp_{\mathrm{in}}$.

We know that $\Pp(D_1(R,T))>0$, and therefore there are numbers $l\in
\N$ and $\Delta,M>0$ such that
%
\begin{equation}
\Pp\bigl(D_1(R,T,l,M,\Delta)\bigr)>0, \label
{eqpositiveprobabilitywithnumberofpoints}
\end{equation}
where
\[
D_1(R,T,l,M,\Delta)=D_1(R,T)\cap\bigl\{N\bigl([0,T]
\times\R\bigr)=N\bigl([\Delta,T-\Delta]\times B_M\bigr)=l\bigr\}.
\]

We let $n\in\N$ be a large number and $\delta\in(0,1/2)$ a small
number to be chosen later and define $L_k=J_k\times
B_\delta$, $k=1,\ldots, n$, where
$J_k=[(2k-1)T/(2n+1), 2kT/(2n+1)]$. The measure $\Pp'_{\mathrm{in}}$ on
configurations in $[0,T]\times\R$ is defined as follows:
all configurations
consist of exactly $n$ independent points, $k$th point distributed
independently in $L_k$ according to
$(2n+1)/(m(B_\delta)T)m(dx)\,dt$. Equivalently,
we can say that $\Pp'_{\mathrm{in}}$ is the distribution of the original
Poissonian point field conditined on having exactly one\vadjust{\goodbreak}
point in each $L_k$.
Thus, the absolute continuity property $\Pp'_{\mathrm{in}}\ll\Pp_{\mathrm{in}}$ holds
and it remains to
prove~\eqref{eqp-prime-positive}. Taking into account
\eqref{eqpositiveprobabilitywithnumberofpoints}, it is
sufficient to show that for any $\omega\in
D_1(R,T,l,M,\Delta)$,
%
\begin{equation}
\Pp'_{\mathrm{in}} \bigl\{\omega'_{\mathrm{in}}
\dvtx\bigl(\omega'_{\mathrm{in}},\omega_{\mathrm{out}}\bigr)\in
D_1(R,T)\cap D_2(R,T) \bigr\}=1. \label{eqwithprob1newomegaisgood}
\end{equation}

First, let us prove that resampled point configurations belong to
$D_1(R,T)$. Since resampling happens only inside
$[0,T]\times\R$,
the time $\tau^-_R$ (depending only on the realization in $(-\infty
,0]\times\R$) does not change.
Therefore $r^-(\omega')<R$ and $\tau^-_R(\omega')<T$, where $\omega
'=(\omega'_{\mathrm{in}},\omega_{\mathrm{out}})$.

Let us prove that $r^+(\omega')<R$ and $\tau^+_R(\omega')<T$. We
need to show that for any $y\in B_R$ and any $t>T$, any
$\gamma'\in M_{0,\R}^{t,y}(\omega')$
satisfies $\gamma'(0)\in B_R$. This is certainly true if $\gamma'$
passes through a point of $\omega'_{\mathrm{in}}$. So, let us
assume that it does not pass through
any points of $\omega'_{\mathrm{in}}$. Therefore, between~$0$ and~$T$ it is a
straight line. Consider now a path $\gamma\in
M_{0,\R}^{t,y}(\omega)$.
We know that $\gamma(0)\in B_R$. Let $t_0=\sup\{t\in[0,T]\dvtx
\gamma
(t)\in B_R\}$. Due to the definition of
$D_1(R,T,l,M,\Delta)$, $t_0>\Delta$.

Let the path $\bar\gamma$ visit all available points in $\omega
'_{\mathrm{in}}$ between $0$ and $t_0/2$, then move straight to
$(t_0,\gamma'(t_0))$
and coincide with $\gamma$ after $t_0$. We are going to show that $\Ac
_{\omega'}^{0,t}(\bar\gamma)<
\Ac_{\omega'}^{0,t}(\gamma')$ so that
$\gamma'$ cannot be a minimizer. Since $\gamma'$ does not pass
through any points of $\omega'_{\mathrm{in}}$,
%
\begin{equation}
\label{eqinterchangingpaths} \Ac_{\omega'}^{0,t}\bigl(
\gamma'\bigr)\ge\Ac_{\omega}^{0,t}\bigl(
\gamma'\bigr) \ge\Ac_{\omega}^{0,t}(\gamma)\ge
\Ac_{\omega'}^{0,t}(\bar\gamma)+\bigl(\Ac_{\omega}^{0,t}(
\gamma)-\Ac_{\omega'}^{0,t}(\bar\gamma)\bigr).
\end{equation}

Let us estimate the difference in the right-hand side. Switching from
$\gamma$ to~$\bar\gamma$, we lose at most $l$ Poissonian
points, but what do we gain?
The action of a path visiting~$r$ points from $\omega'_{\mathrm{in}}$ in a row
does not exceed
\[
A(r)=r\frac{(2\delta)^2}{2T/(2n+1)}-r,
\]
and, since there are at least $nt_0/(3T)$ points visited by $\bar
\gamma$ in $\omega'_{\mathrm{in}}$ between $0$ and $t_0$,
\[
\Ac_{\omega'}^{0,t_0}(\bar\gamma)<\frac{nt_0}{3T} \biggl(
\frac
{(2\delta)^2}{2T/(2n+1)}-1 \biggr)+\frac{R^2}{2(t_0/2)}.
\]
Therefore,
\[
\Ac_{\omega}^{0,t}(\gamma)-\Ac_{\omega'}^{0,t}(
\bar\gamma)\ge-l + \frac{nt_0}{3T} \biggl(1-\frac{(2\delta
)^2}{2T/(2n+1)} \biggr)-
\frac
{R^2}{2(t_0/2)}.
\]
Choosing $n$ to be large and $\delta$ small, we see that the
right-hand side is positive, which in conjunction
with~\eqref{eqinterchangingpaths}, gives
the desired inequality $\Ac_{\omega'}^{0,t}(\bar\gamma)< \Ac
_{\omega'}^{0,t}(\gamma')$.

This finishes the proof of $r^+(\omega')<R$ and $\tau^+_R(\omega
')<T$. It is also easy to adjust the above argument to
show that
$r(\omega')<R$ and $\tau_R(\omega')<T$, and in the same way one can
prove that $r(\theta^T\omega')<R$ and
$\tau_R(\theta^T\omega')<T$. Thus, $\omega'\in D_1(R,T)$, and
it remains to prove that $\omega'\in D_2(R,T)$ a.s.\vadjust{\goodbreak}

Let us prove the following claim: for any points $x,y\in B_R$, an
optimal path $\gamma\in M_{x,0}^{y,T}(\omega')$ cannot
avoid all
points of $\omega'$ between $0$ and $T/3$. In fact, if it does avoid
these points, then
%
\begin{equation}
\label{eqactionforonethirdavoiding} \Ac_{\omega'}^{0,T}(\gamma)
\ge-\tfrac{2}{3}n-1.
\end{equation}
On the other hand, consider the path $\bar\gamma\in\Gamma
_{0,x}^{T,y}$ that visits all points of $\omega'$ between
$T/8$ and $7T/8$,
%
\begin{equation}
\label{eqactionforthreequartersvisits} \Ac_{\omega'}^{0,T}(\bar
\gamma)\le2 \frac{(2R)^2}{2
T/8}- \biggl(n \biggl(\frac{7T}{8}-
\frac{T}{8} \biggr)-2 \biggr) \biggl(1-\frac{(2\delta
)^2}{2T/(2n+1)} \biggr),
\end{equation}
and, for sufficiently small $\delta$ and large $n$, $\Ac_{\omega
'}^{0,T}(\bar\gamma)<\Ac_{\omega'}^{0,T}(\gamma)$, which
contradicts our assumption
$\gamma\in M_{x,0}^{y,T}(\omega')$. Our claim is proven, and in the
same way one can prove that an optimal path
$\gamma\in M_{x,0}^{y,T}(\omega')$ must pass
through one of the points of $\omega'$ between $2T/3$ and $T$.
Clearly, for sufficiently small $\delta>0$ any optimal
path passing through
a point in $\omega'$ between $0$ and $T/3$ and a point in $\omega'$
between $2T/3$ and $T$ also passes through all
points in $\omega'$ in between. Therefore, $\omega'\in D_2(R,T)$, and
the proof of Lemma~\ref{lmcollapse-with-positive-prob} is
complete.
\end{pf}


We can now construct a global solution. For a set $A$ of paths and a
time interval $[s,t]$, we denote by
$A |_{[s,t]}$
the set of restrictions of all trajectories from $A$ to $[s,t]$.

\begin{lemma} Let $R>0$. There are two random times $\sigma_0,\sigma
_1>0$ such that for all $x\in B_R$ and any two times
$t_1,t_2>\sigma_1$,
\[
M_{-t_1,\R}^{0,x}(\omega) |_{[-\sigma_0,0]}= M_{-t_2,\R
}^{0,x}(
\omega) |_{[-\sigma_0,0]}.
\]
\end{lemma}

\begin{pf}The ergodicity of $\theta^1$,
Lemma~\ref{lmcollapse-with-positive-prob} and Poincar\'e's Recurrence
theorem imply that with
probability $1$, there is an integer time $n>T$ such that $\theta
^{-n}\omega\in D_1(R,T)\cap D_2(R,T)$.
Without loss of generality we can assume that $R>R'$, where $R'$ is defined
in Lemma~\ref{lmmain-loc}. The last part of that lemma implies that
if we define
%
\begin{equation}
\sigma_1(\omega)=n+\bar\tau_R\bigl(
\theta^{-n}\omega\bigr)+\bar\tau_R\bigl(
\theta^{-n+T}\omega\bigr), \label{eqcouplingtime}
\end{equation}
then for any $t>\sigma_1$ and any $x\in B_R$, any path $\gamma\in
M_{-t,\R}^{0,x}(\omega)$
satisfies $\gamma(-n)\in B_{r^{\pm}(\theta^{-n}\omega)}\subset B_R
$ and $\gamma(-n+T)\in
B_{r^{\pm}(\theta^{-n+T}\omega)}\subset B_R$.
Here we used the fact that $\theta^{-n}\omega\in D_1(R,T)$.

Since $\theta^{-n}\omega\in D_2(R,T)$, for any $t_1$ and $t_2$
satisfying~\eqref{eqcouplingtime} and any $x\in B_R$,
any paths $\gamma_1\in M_{-t_1,\R}^{0,x}(\omega)$
and $\gamma_2\in M_{-t_2,\R}^{0,x}(\omega)$
pass through a common point at some time $\sigma_0\in[-n,-n+T]$.
Therefore, the restrictions of the sets of minimizers
on $[-\sigma_0,0]$ coincide, and the proof is complete.
\end{pf}

\begin{remark}\label{remone-sidedminimizers}
In fact, there is an infinite, strictly increasing sequence
$(n_k)_{k\in\N}$ such that
$\theta^{-n_k}\omega\in D_1(R,T)\cap D_2(R,T)$ for all $k$.
Therefore, the theorem can be strengthened. Its conclusion
holds for any (random
or deterministic) $\sigma_0>0$. In particular, the finite time
minimizers stabilize to a limiting infinite one-sided
minimizer.
\end{remark}

We can now finish our construction of the global solution $u$. For any
$r\in\R$ and $R>0$, the restriction of
$\Phi^{-s,t}0$ on $B_R$ stabilizes for large values of $s$,
and
\[
u_\omega(t)=(\Ub,d)\lim_{s\to\infty} \Phi^{-s,t}0
\]
is well defined. Clearly, the construction of $u_\omega(t)$ depends
only on
the restriction of $\omega$ on $(-\infty,t]\times\R$.

Since restrictions of minimizers are also minimizers, we can deduce
that for any time interval $[t_0,t_1]$,
$u_\omega(t_1)$ is the solution of the Cauchy problem with initial
value $u_\omega(t_0)$. Therefore, thus constructed
function $u_\omega$
is a global solution of the Burgers equation corresponding to the
realization of the random forcing~$\omega$. This
proves parts (\hyperlink{itmeasurability-of-u}{1}) and
(\hyperlink{itglobal-solution}{2})
of the main theorem.

Let us prove part~(\hyperlink{itpiecewise-linear}{3}) of the main theorem. For
each point $x$ of continuity of~$u_\omega$,
we denote $\pi_\omega(x)$ the Poissonian point that is visited last
by the
minimizer $\gamma\in M_{-t,\R}^{0,x}$ for sufficiently large $t$. The
map $\pi_\omega$ is piecewise constant.
If $\pi_\omega(x)=(s_i,x_i)$ for all $x$ in an interval $J$, then
\[
u_\omega(x)=\frac{x-x_i}{|s_i|},\qquad x\in J,
\]
and part~(\hyperlink{itpiecewise-linear}{3}) follows.

\section{\texorpdfstring{The behavior of global solution $u(t,x)$ as $x\to\infty$}
{The behavior of global solution u(t,x) as x to infinity}}
\label{secatinfty}
In this section we prove part~(\hyperlink{itrho}{4}) of the main theorem. We
will concentrate on proving the limit behavior as
$x\to+\infty$ since
the limit $x\to-\infty$ can be studied in exactly the same way.

The idea is that if we want to consider, say, a path in $M_{0,\R
}^{T,x}$ for large values of $x$ and $T$, then the path
naturally
decomposes into two parts.
Most Poissonian points are scattered over a compact domain, so
in a certain time interval $[0,t]$ the path mostly stays in a compact
domain around the origin collecting action at
approximately linear rate $S<0$,
and then it leaps from the compact domain straight to $x$ roughly with
constant speed between $t$ and $T$, hardly
meeting any Poissonian points in this regime and
collecting approximately $x^2/(2(T-t))$ action. Finding the minimum of
\[
St+\frac{x^2}{2(T-t)},\qquad t\in[0,T],
\]
we obtain that the optimal $t$ satisfies
\[
\frac{x}{t}=\sqrt{-2S}.\vadjust{\goodbreak}
\]

This nonrigorous argument shows that we can hope that part~(\hyperlink{itrho}{4}) of the main theorem holds with $q=\sqrt{-2S}$.

To make this argument precise, we need to control the deviations of the
action on $[0,t]$ from $St$ and to control the
behavior of minimizers very far from
the origin where Poissonian points are sparse.

For any connected set $I\subset\R$ and any time interval $[t_0,t_1]$,
we can define

\[
S_I^{t_0,t_1}=S_I^{t_0,t_1}(\omega)=\inf
\bigl\{\Ac^{t_0,t_1}_\omega(\gamma)\dvtx\gamma(s)\in I\mbox{ for
all } s\in[t_0,t_1] \bigr\}.
\]
Clearly, this function is superadditive:
\[
S_I^{t_0,t_2}\ge S_I^{t_0,t_1}+
S_I^{t_1,t_2},\qquad t_0\le t_1 \le
t_2.
\]

The ergodicity of the flow $(\theta^t)_{t\in\R}$ and Kingman's
subbaditive ergodic theorem imply that the following
random variable is well defined and a.s.-constant:
\[
S_I=\lim_{t_1\to\infty}\frac{1}{t_1-t_0}S_I^{t_0,t_1}=
\lim_{t_0\to-\infty}\frac{1}{t_1-t_0}S_I^{t_0,t_1}=\mathop{
\lim_{t_0\to-\infty}}_
{t_1\to\infty}\frac{1}{t_1-t_0}S_I^{t_0,t_1}.
\]
Clearly, $S_{B_R}$ is a nonincreasing negative function of $R>0$, and
we define $S=\lim_{R\to\infty}S_{B_R}<0$.

\begin{lemma}\label{lmSaslimit} Thus defined constant $S$ satisfies
$S=S_\R$.
\end{lemma}
\begin{pf}
Obviously, $S_\R\le S_{B_R}$ for any $R$.
Therefore, we only have
to prove that $S_\R\ge S$.
Let us take any $t>0$ and any path $\gamma$ realizing $S_{\R}^{0,t}$.
Taking any $R>0$ and decomposing $\gamma$ into
parts that stay inside $B_R$ and outside $B_R$, we
see that
\[
S_\R^{0,t}=\Ac^{0,t}_\omega(\gamma)\ge
S_{B_R}^{0,t}-N\bigl([0,t]\times B_R^c
\bigr).
\]
Dividing by $t$ and taking $t\to\infty$, we obtain
\[
S_\R\ge S_{B_R}-m\bigl(B_R^c
\bigr).
\]
Taking $R\to\infty$ finishes the proof of the lemma.
\end{pf}

\begin{lemma}\label{lmexpectation-converge}
\[
\lim_{t\to\infty} \frac{\E S_{B_R}^{0,t}}{t} = S_{B_R}.
\]
\end{lemma}
\begin{pf}Since
\[
0\ge\frac{S_{B_R}^{0,t}}{t} \ge-\frac{N([0,t]\times\R)}{t},
\]
the lemma follows by dominated convergence.
\end{pf}

\begin{lemma} \label{lmlarge-deviation}
\begin{longlist}[(1)]
\item[(1)] Let $R\in(0,\infty]$. If $S'<S$, then there are
constants $c=c(S')>0$ and $T_0=T_0(S')>0$ such
that for $T>T_0$,
\[
\Pp\biggl\{\inf_{t\ge T} \frac{S^{0,t}_{B_R}}{t}<S' \biggr
\}<e^{-cT}.
\]
\item[(2)] If $S'>S$, then there is a constant $R_0=R_0(S')$ with the
following property: for every $R\in[R_0,\infty]$, there
are $c=c(S', R)>0$ and $T_0=T_0(S',\break R)>0$ such that for $T>T_0$,
\[
\Pp\biggl\{\sup_{t\ge T} \frac{S_{B_R}^{0,t}}{t}> S' \biggr
\}<e^{-cT}.
\]
%
\end{longlist}
\end{lemma}

\begin{pf}Recall that $S_{B_R}\downarrow S<0$.
Since $S'<S<S_{B_R}$, Lemma~\ref{lmexpectation-converge} allows us to
choose $s$ such that
%
\begin{equation}
\label{eqestimateonexpect} \E S_{B_R}^{0,s}> s
\frac{S_{B_R}+S'}{2}.
\end{equation}
Then we notice that
\begin{eqnarray*}
S_{B_R}^{0,t}&\ge& S_{B_R}^{0,s}(
\omega)+S_{B_R}^{s,2s}(\omega)+\cdots+S_{B_R}^{[t/s]s,([t/s]+1)s}(
\omega)
\\
&\ge& S_{B_R}^{0,s}(\omega)+S_{B_R}^{0,s}
\bigl(\theta^{s}\omega\bigr)+\cdots+S_{B_R}^{0,s}
\bigl(\theta^{[t/s]s}\omega\bigr).
\end{eqnarray*}
Let us denote the right-hand side by $\Sigma_{[t/s]+1}$. It is the
sum of $[t/s]+1$ i.i.d. nonpositive random variables with
finite exponential moments, and with expectations estimated by
\eqref{eqestimateonexpect}.

Since $\frac{S't}{[t/s]+1}\to S's$ as $t\to\infty$, and
$S'<(S_{B_R}+S')/2$, the estimate
\[
\Pp\bigl\{S_{B_R}^{0,t}<S't\bigr\}\le
K_1e^{-c_1t}
\]
for all $S'<S$, some $K_1=K_1(S')$, $c_1=c_1(S')>0$ and all $t>0$,
is a consequence of the classical Cram\'er large
deviation estimate.
Since $S_{B_R}^{0,t}$ is nonincreasing in $t$, and
\[
\inf_{s\in[t,t+1]}\bigl(S_{B_R}^{0,s}-S_{B_R}^{0,t}
\bigr)\ge- N\bigl([t,t+1]\times\R\bigr),
\]
we can use this maximal inequality in a standard way to interpolate
between~$t$ and $t+1$ and obtain
\[
\Pp\biggl\{\inf_{s\in[t,t+1]} \frac{S_{B_R}^{0,s}}{s}<S' \biggr\}
\le
K_2e^{-c_2t}
\]
for all $S'<S$, some $K_2=K_2(S')$, $c_2=c_2(S')>0$ and all $t>0$. Now
the first part of the lemma follows.
The proof of the second part of the lemma is essentially the same.
\end{pf}

Now we turn to ruling out paths ending at a large $x$ and having slopes
deviating significantly from $q$.
For any $x>0$, $\eps\in(0, q)$, and $R\in(0,x)$ let us denote
\begin{eqnarray*}
Q_+(x,\eps,R)&=&\bigl\{(t,y)\dvtx y\in(R,x), t<0, y< x+t(q+\eps)\bigr\},
\\
Q_-(x,\eps,R)&=&\bigl\{(t,y)\dvtx y\in(R,x), t<0, y> x+t(q-\eps)\bigr\}.
\end{eqnarray*}

\begin{lemma} \label{lmnoshortpath}
For each $\eps>0$, there is random variable $R$ such that:
\begin{longlist}[(1)]
\item[(1)] \label{lmQplus}For any
$x>2R$, a path $\gamma$ inside $Q_+$ connecting a point $(t,y)\in
\partial Q_+(x,\eps,  R)$ with $y> R$ to $(0,x)$ and
$\dot\gamma(0)>q+2\eps$
cannot be a minimizer.
\item[(2)] For any
$x>2R$, a path $\gamma$ inside $Q_-$ connecting a point $(t,y)\in
\partial Q_-(x,\eps,  R)$ with $y> R$ to $(0,x)$ and
$\dot\gamma(0)<q-2\eps$
cannot be a minimizer.
\end{longlist}
\end{lemma}

\begin{pf}First, we notice that condition \eqref
{eqmomentcondition} implies
that the set $\{(t,x)\dvtx\break -x<t<0 \}$
contains finitely many Poissonian points with probability $1$.
Therefore, we can define a random variable $r_0$ such
that
with probability $1$, there are no Poissonian points in the set
\[
A=\bigl\{(t,x)\dvtx x>r_0, -x<t<0 \bigr\}.
\]
Let us now define $\rho=q/(2+q)$ and, for any $x$,
\[
A_x=[-\rho x/q,0]\times[x-\rho x,x].
\]
It is easy to check that $(-\rho x/q, x-\rho x)\in A$ for sufficiently
large values of~$x$. Therefore, for these values
of $x$, $A_x\subset A$.
We conclude that there is a random number $r_1$ such that for $x>r_1$,
there are no Poissonian
points in~$A_x$.

Let us now take a path $\gamma$ satisfying the conditions of part~\ref
{lmQplus} of the lemma.
We would like to compare this path to the straight line segment
connecting $(t,y)$ and $(0,x)$.

\begin{lemma} \label{lmdeviation-from-straight} Consider points
$(t_0,x_0)$, $(t_1,x_1)$, $(t_2,x_2)$ satisfying
$t_0<t_1<t_2$.
The free action (i.e., the action without taking account the
contribution from the Poissonian points) of the path
connecting these points is bounded below by
\[
\frac{(x_2-x_0)^2}{2(t_2-t_0)}+\frac{2}{t_2-t_0}(x_1-\bar x)^2,
\]
where
\[
\bar x=\frac{x_2(t_1-t_0)+x_0(t_2-t_1)}{t_2-t_0},
\]
so that $|x_1-\bar x|$ is the distance from $(t_1,x_1)$ to the straight
line connecting $(t_0,x_0)$ and $(t_2,x_2)$,
measured along the $x$-axis.
\end{lemma}

\begin{pf}The free action minimizing path consists of
two straight line
segments connecting $(t_0,x_0)$ to $(t_1,x_1)$ and
$(t_1,x_1)$ to $(t_2,x_2)$. The resulting\vadjust{\goodbreak} action is
a quadratic polynomial in $x_1$,
\begin{eqnarray*}
f(x_1)&=&\frac{(x_1-x_0)^2}{2(t_1-t_0)}+\frac
{(x_2-x_1)^2}{2(t_2-t_1)}
\\
&=&\frac{t_2-t_0}{2(t_2-t_1)(t_1-t_0)}(x_1-\bar x)^2+\frac
{(x_2-x_0)^2}{2(t_2-t_0)},
\end{eqnarray*}
and the estimate
\[
\frac{t_2-t_0}{2(t_2-t_1)(t_1-t_0)}=\frac
{(t_2-t_1)+(t_1-t_0)}{2(t_2-t_1)(t_1-t_0)}\ge\frac
{2}{(t_2-t_1)+(t_1-t_0)} =\frac{2}{t_2-t_0}
\]
completes the proof.
\end{pf}

Let us denote by $(s,z)$ the Poissonian point that is connected by the
last segment of path $\gamma$ to $(0,x)$.
To apply Lemma~\ref{lmdeviation-from-straight},
we must estimate the distance from $(s,z)$ to the straignt line
connecting $(t,y)$ and $(0,x)$ measured along the
$x$-axis, that is, $x-z+s(q+\eps)$.
Since there are no Poissonian points in $A_x$,
we have $(s,z)\in Q_+\setminus A_x$ and, consequently, $z\le x-\rho x$.
Since $\dot\gamma(0)>q+2\eps$, we have
$(x-z)/(-s)>q+2\eps$. The minimum of $x-z+s(q+\eps)$
under these restrictions is attained at $(-\rho x/(q+2\eps), x-\rho
x)$ and equals $\eps\rho x/(q +2 \eps)$. We also
have $|t|<(x-R)/(q +\eps)$.
Therefore, Lemma~\ref{lmdeviation-from-straight} implies that the
action gain of $\gamma$ compared to the straight line
motion between $t$, and $0$
is at least
%
\begin{equation}
\label{eqlowerboundonactionincrement} \frac{2(q +\eps)}{x-R}\frac
{\eps^2\rho^2 x^2}{(q + 2\eps)^2}\ge K(
\eps)x,\qquad x>(2R)\vee r_1(\omega)
\end{equation}
for some $K(\eps)$.

Now we must estimate the effect of Poissonian points.
We use Lemma~\ref{lmlln-bound} to find $R_0=R_0(\omega)$ such that
$N([-t,0]\times B_{R_0}^c)< t K(\eps)q/2$ for all $t>0$. Since the
time component of any point in $Q_+$ is bounded by
$x/q$ in absolute value, we see that if $R>R_0$, then there are at most
$(x/q) K(\eps)q/2< K(\eps)x/2$ Poissonian points
in $Q_+$.
Therefore, the reduction of action due to visits to Poissonian points
does not exceed $K(\eps)x/2$, and cannot
compensate for the action gain computed
in~\eqref{eqlowerboundonactionincrement}. Therefore, if we choose
$R>R_0\vee r_1$, then for any $x>2R$ the straight
line segment from $(t,y)$ to $(0,x)$ is more efficient than any path
$\gamma$ satisfying
the imposed requirements, so $\gamma$ cannot be a minimizer. The proof
of the first part of the lemma is complete.
The proof of the second part is similar, and we omit it.
\end{pf}

\begin{lemma}\label{lmexclude-lower-paths-1} For any $\eps>0$ there
are random variables $R>0$ and $X>0$ such that if
$x>X$ and $\tau>(x-R)/(q-\eps)$, then
no path $\gamma\in M_{-\tau,R}^{0,x}$ can satisfy
%
\begin{equation}
\label{eqgammastaysoutside} \gamma(s)> \bigl(x+(q-\eps)s \bigr
)\vee R, \qquad s\in(-
\tau,0 ].
\end{equation}
\end{lemma}

\begin{pf}We begin with taking $\delta>0$ (to be chosen later)
and using Lemma~\ref{lmlln-bound} to find $R_0$ such that for all
$\tau>(x-R)/(q-\eps)$ and all $R>R_0$, the action of
any path $\gamma$ satisfying~\eqref{eqgammastaysoutside},
connecting $(-\tau,R)$ to $(0,x)$ and staying outside of $B_R$ for all
times in $(-\tau,0]$,
is at least
\[
\frac{(x-R)^2}{2\tau}-\delta\tau.
\]
To prove that $\gamma\notin M_{-\tau,R}^{0,x}$, let us find a better
path $\tilde\gamma$ in $\Gamma_{-\tau,R}^{0,x}$.
First, we will choose $\tilde\gamma$ so that $\tilde\gamma|_{[-\tau
+1,-[(x-R)/q]-2}\in
M_{-\tau+1,B_R}^{-[(x-R)/q]-2,B_R}$. Then we denote $x_1=\tilde\gamma
(-[(x-R)/q]-2)$ and $x_2=\tilde\gamma(-\tau+1)$.
The remaining parts of $\tilde\gamma$ are straight line segments
connecting $(-\tau,R)$ to $(-\tau+1,x_2)$,
$(-[(x-R)/q]-2,x_1)$ to $(-[(x-R)/q]-1,R)$, and $(-[(x-R)/q]-1,R)$ to $(0,x)$.

The action of this path is at most
\[
\frac{(x-R)^2}{2([(x-R)/q]+1)}+\frac{(2R)^2}{2}+ S_{B_R}^{0,\tau
-1-([(x-R)/q]+2)}\bigl(
\theta^{-[(x-R)/q]-2}\omega\bigr)+\frac
{(2R)^2}{2}.
\]
We want to exclude the situation where
%
\begin{equation}
\Ac^{-\tau,0}(\tilde\gamma)\ge\Ac^{-\tau,0}(\gamma).
\label{eqineq-wrong-sign}
\end{equation}
Suppose~\eqref{eqineq-wrong-sign} holds. Then
\[
\frac{(x-R)^2}{2([(x-R)/q]+1)}+(2R)^2+ S_\tau\ge\frac
{(x-R)^2}{2\tau}-
\delta\tau,
\]
where we denoted $S_\tau=S_{B_R}^{0,\tau-3-[(x-R)/q]}(\theta
^{-[(x-R)/q]-2}\omega)$ for brevity.
This can be rewritten as
%
\begin{equation}
\label{eqneed-large-dev} \frac{S_\tau}{\tau-3-[(x-R)/q]}\ge U,
\end{equation}
where
\[
U=\frac{({(x-R)^2}/{2}) ({1}/{\tau}-
{1}/{([(x-R)/q]+1)} )-\delta\tau-(2R)^2}{\tau-3-[(x-R)/q]}.
\]
To apply large deviation estimates from Lemma~\ref
{lmlarge-deviation}, we need to estimate $U$ and the length of time
interval in the definition
of $S_\tau$. Since $\gamma$ satisfies \eqref{eqgammastaysoutside}
for $s\in[-\tau,0]$, we have
$\tau>(x-R)/(q-\eps)$.
For sufficiently small $\eps$,
%
\begin{eqnarray}
\label{eqlowerboundontau} \tau-3-\bigl[(x-R)/q\bigr]&\ge&\frac
{x-R}{q-\eps}-
\frac{x-R}{q}-4 \ge\frac
{(x-R)\eps}{q(q-\eps)}-4
\nonumber
\\[-8pt]
\\[-8pt]
\nonumber
&\ge&\frac{(x-R)\eps}{2q^2}.
\end{eqnarray}
To estimate $U$, we first notice that, due to \eqref
{eqlowerboundontau}, there is a random variable
$X_1(R,\eps,\delta)$ such that $x>X_1$ implies
%
\begin{equation}
\label{eqauxil-contrib} \frac{\delta\tau+(2R)^2}{\tau
-3-[(x-R)/q]}<2\delta.
\end{equation}
For the same reason, there is a random variable $X_2(R,\eps,\delta)$
such that $x>X_2$ implies
%
\begin{eqnarray}\label{eqmain-contrib}
\nonumber
&&\frac{({(x-R)^2}/{2} )({1}/{\tau}-
{1}/{([(x-R)/q]+1)} )}{\tau-3-[(x-R)/q]}\\
&&\qquad=\frac{(x-R)^2}{2\tau([
(x-R)/q]+1)}\cdot\biggl(-\frac{\tau-1-[(x-R)/q]}{\tau
-3-[(x-R)/q]}
\biggr)
\nonumber
\\[-8pt]
\\[-8pt]
\nonumber
&&\qquad\ge-\frac{(x-R)q}{2\tau}(1+\delta)
\\
&&\qquad\ge-\frac{q(q-\eps)}{2}(1+\delta) ,\nonumber
\end{eqnarray}
where the last inequality follows from $\tau>(x-R)/(q-\eps)$.

Combining \eqref{eqauxil-contrib} and~\eqref{eqmain-contrib} and
choosing $\delta$ sufficiently small we see that
(with the choices of $R$ and $x$ described above)
%
\begin{equation}
U>-\frac{q(q-\eps/2)}{2}. \label{eqlowerboundonU}
\end{equation}

Notice that if $k$ is sufficiently large, all the above estimates apply
uniformly for all $x\in[R+kq,R+(k+1)q]$ and all
$\tau\ge(x-R)/(q-\eps)$.

Let us denote by $B_k$ the event that for some $x\in[R+kq,R+(k+1)q]$
and some $\tau\ge(x-R)/(q-\eps)$ there is a path
$\gamma\in M_{-\tau,R}^{0,x}$
satisfying~\eqref{eqgammastaysoutside} for $s\in[-\tau,0]$.
The definition of $q$, inequality~\eqref{eqlowerboundonU}, and
Lemma~\ref{lmlarge-deviation} imply that for some
$c>0$ and all sufficiently large $k$,
\[
\Pp(B_k)<e^{-ck}.
\]
Now the Borel--Cantelli lemma implies that with probability, 1 only
finitely many events $B_k$ happen, and the proof is
complete.
\end{pf}

\begin{lemma}
\label{lmexcludelongerlower}
There are positive random variables $R,X$ and $(T_x)_{x>0}$ such that
if $x>X$ and $\tau>T_x$, then for any $y\in\R$,
no path $\gamma\in M_{-\tau,y}^{0,x}$ can satisfy~\eqref
{eqgammastaysoutside}.
\end{lemma}

\begin{pf}
If for some $y$ a path $\gamma\in\Gamma
_{-\tau,y}^{0,x}$
satisfies~\eqref{eqgammastaysoutside} and the time
$\tau'=\sup\{s\dvtx\gamma(s)\le R\}$ is well defined,
then we can apply Lemma~\ref{lmexclude-lower-paths-1} with $\tau$
replaced by $\tau'$ to see that $\gamma$ cannot be a
minimizer for appropriately
chosen~$R$ and $X$.

Let us fix $\delta\in(0,-S)$. Due to Lemma~\ref{lmlln-bound}, we
can choose $R$ large enough to ensure that
$\Ac^{-\tau,0}>-\delta\tau$ for
any $\gamma$ satisfying $\gamma(s)>R$ for all $s\in[-\tau,0]$. On
the other hand, the optimal
action is asymptotic to $S\tau$ as $\tau\to\infty$, so $\gamma$
cannot be a minimizer for large values of $\tau$.
\end{pf}

\begin{lemma} \label{lmexcludelongerupper}There are positive random
variables $R,X$ and $(T_x)_{x>0}$ such that for
$x>X$ and $T>T_x$,
and any $y\in\R$ no $\gamma\in M_{-T,y}^{0,x}$
can satisfy
%
\begin{equation}
\label{eqabovezone} \gamma(s)<x+s(q+\eps),\qquad  s\in\bigl
[-(x-R)/(q+\eps),0\bigr].
\end{equation}
\end{lemma}
\begin{pf}We need an auxiliary path $\bg\in\Gamma
_{-(1+\eps
)[(x-R)/q],-R}^{0,x}$. This special path consists of three straight
line segments connecting
consecutively $(-(1+\eps)[(x-R)/q],-R)$ to $(-[(x-R)/q],-R)$ to
$(-[(x-R)/q]+1,R)$ to $(0,x)$.
\end{pf}

\begin{lemma} There are positive random variables $R$ and $X$ such that
for any $x>X$, any $T>-(1+\eps)[(x-R)/q]$ and
any $y\in\R$,
any $\gamma\in M_{-T,y}^{0,x}$ satisfying~\eqref{eqabovezone}
intersects $\bg$.
\end{lemma}

\begin{pf}Denote by $B_k$ the event that there are $y\in
\R$ and $\gamma
\in M_{-T,y}^{-t,-R}$ for some $t<k$ and $T>(1+\eps)k$
such that
$\gamma(s)<-R$ for all $s\in[-T,-t]$.

Lemmas~\ref{lmlln-bound} and~\ref{lmlarge-deviation} imply that for
sufficiently large $R$, for
some constants $c_1,c_2>0$ and all $k$, $P(B_k)\le c_1e^{-c_2\eps k}$.
The Borel--Cantelli lemma implies that with
probability 1, only finitely many events $B_k$
happen.

Clearly, if there is a path $\gamma$ satisfying the conditions of the
lemma and not intersecting $\bg$, then
$B_k$ holds for $k=[(x-R)/q]$. Since only finitely many $B_k$ can hold,
a path with these properties is impossible for
sufficiently large~$x$.
\end{pf}

\begin{lemma} There are random variables $R$ and $X$ such that for
$x>X$ and any path $\gamma$
satisfying~\eqref{eqabovezone}
and intersecting $\bg$ at some time $-\tau$,
%
\begin{equation}
\Ac^{-\tau,0}(\gamma)>\Ac^{-\tau,0}(\bg). \label{eqcanimproveaction}
\end{equation}
\end{lemma}
\begin{pf}
First let us consider the possibility that
$\tau<[(x-R)/q]$. We denote
$\nu=\inf\{s\dvtx\gamma(-s)=R\}$. For any $\delta>0$,
there is $R$ such that for $x>R$,
%
\begin{equation}
\label{eqlowerestimateonAgamma} \Ac^{-\tau,0}(\gamma)\ge
\frac{(x-R)^2}{2\nu}+S_{\R}^{-[
{(x-R)}/{q}]+1,-\nu}-\delta\nu.
\end{equation}

On the other hand,
\[
\Ac^{-\tau,0}(\bg)\le\frac{(x-R)^2}{2([{(x-R)}/{q}]-1)}+\frac
{(2R)^2}{2}.
\]
If \eqref{eqcanimproveaction} is violated, the last two
inequalities imply
%
\begin{eqnarray}\label{eqimplicationofwronineq}
\qquad\frac{S_{\R}^{-[{(x-R)}/{q}]+1,-\nu}}{[{(x-R)}/{q}]-1-\nu
}&<&-\frac{(x-R)^2}{2\nu([({x-R)}/{q}]-1)}+\frac{\delta\nu
}{{[{(x-R)}/{q}]-1-\nu}}
\nonumber
\\[-8pt]
\\[-8pt]
\nonumber
&&{}+ \frac{2R^2}{{[{(x-R)}/{q}]-1-\nu}}.
\end{eqnarray}
From \eqref{eqabovezone} we know that $\nu<(x-R)/(q+\eps)$. We can
use this to derive that the second term in the
right-hand side is bounded by $K\delta/\eps$
for a constant $K>0$ and the third term converges to $0$ as $x\to
\infty$.
Choosing $\delta$ sufficiently small, then choosing $R$ so that
\eqref{eqlowerestimateonAgamma} holds, we conclude that for
sufficiently large $x$, the right-hand side does not exceed
$-q(q+\eps/2)/2$.
Now the large deviation estimate of Lemma~\ref{lmlarge-deviation} and
the Borel--Cantelli lemma imply that
\eqref{eqimplicationofwronineq}
can hold true only for a bounded set of $x$.

Now we have to exclude the paths $\gamma$ that cross $\bg$ for the
first time at~$-R$.
By considering a smaller value of $\delta$ in the above reasoning, it
is easy to strengthen it
and conclude that that there is $\Delta>0$ such that for sufficiently
large $x$, all paths $\gamma$ satisfying
this restriction satisfy also
%
\begin{equation}
\Ac^{-[(x-R)/q],0}(\gamma)>\Ac^{-[(x-R)/q],0}(\bg)+\Delta(x-R).
\label{eqDelta}
\end{equation}
On the other hand, denoting $I_{\eps,R,x}= [-(1+\eps)[\frac
{x-R}{q}],-[\frac{x-R}{q}] ]$,
\[
\inf_{t\in I_{\eps,R,x} }\Ac^{t,[(x-R)/q]}(\gamma)\ge- N\bigl(I_{\eps
,R,x}\times(-
\infty,-R]\bigr).
\]
Suppose $R$ is chosen so that $\E N(I_{\eps,R,x}\times(-\infty
,-R])<\Delta(x-R)/2$. Then
probability that for some $x\in[k,k+1]$, there is a path $\gamma$,
with $\gamma(0)=x$, satisfying~\eqref{eqDelta}
and violating \eqref{eqcanimproveaction},
decays exponentially in $k$. An application of the Borel--Cantelli
finishes the proof.
\end{pf}

Part~(\hyperlink{itrho}{4}) of the main theorem follows now from Lemmas~\ref
{lmnoshortpath},~\ref{lmexclude-lower-paths-1},
\ref{lmexcludelongerlower}
and~\ref{lmexcludelongerupper}.

\section{The global solution as a one-pont attractor}\label{secattractor}
In this section we prove part~(\hyperlink{itattractor}{5}) of the main theorem.

Let us denote by $M_{t_0,\R,V}^{t_1,x}$ the set of minimizers of
\eqref{eqminimizationproblem}.

\begin{lemma}\label{lminitialpointloc} Suppose $V$ satisfies~\eqref
{eqasymptoticslope}. Then for any $L>0$,
there is a random variable $R_0>0$ such that for all $t>0$ and all
$x\in B_L$, any $\gamma\in M_{0,\R,V}^{t,x}$
satisfies $\gamma(0)\in B_{R_0}$.
\end{lemma}

\begin{pf}Property~\eqref{eqasymptoticslope} implies
that there is
$\alpha\in(0,q)$ such that
%
\begin{equation}
\label{eqVforlargey} V(y)> -\alpha y \qquad\mbox{for sufficiently large } y>0.
\end{equation}

Let us take a small $\delta>0$ to be chosen precisely later and use
Lemmas~\ref{lmSaslimit} and~\ref{lmlln-bound} to
find $h>L$ such that
$S_{B_h}<S+\delta$ and $N(B_h^c\times[0,s])<\delta s$ for all $s>0$.
Let us consider $y\in[k,k+1]$ for large values of
$k\in\N$ and
estimate the action of a path $\gamma\in\Gamma_{0,y}^{t,x}$. Since
$x\in B_L\subset B_h$, we can define
\[
\tau=\inf\bigl\{s\dvtx\gamma(s)\in B_h\bigr\}.
\]
The complete action of this path on $[0,\tau]$ satisfies
\[
\Ac_V^{0,\tau}(\gamma)\ge-\alpha(k+1)-\delta\tau+
\frac
{(k-h)^2}{2\tau}.
\]
On the other hand, there is a number $C(h)$ such that the optimal path
$\gamma_h\in M_{0,B_h}^{\tau,h}$
satisfies
\[
\Ac_V^{0,\tau}(\gamma_h)\le C(h)+S_{B_h}^{0,\tau-1}.
\]
Therefore, if $\gamma$ is optimal, then
\[
-\alpha(k+1)-\delta\tau+\frac{(k-h)^2}{2\tau}\le
C(h)+S_{B_h}^{0,\tau-1}.
\]
According to the definition of $S_{B_h}$, there is $\tau_{h,\delta}$
such that if $\tau>\tau_{h,\delta}$,
then $S_{B_h}^{0,\tau-1}+C(h)\le(S_{B_h}+\delta)\tau\le(S+2\delta
)\tau$. Therefore, if for optimal $\gamma$,
$\tau>\tau_{h,\delta}$,
then
%
\begin{equation}
(S+3\delta)\tau- \frac{(k-h)^2}{2\tau} \ge-\alpha(k+1). \label
{eqfindoptimalpathtomiddle}
\end{equation}

Elementary calculus shows that the global maximum of the left-hand
side in $\tau$ is achieved at
$\tau^*=(k-h)/\sqrt{2(-S-3\delta)}$
and equals $-(k-h)\sqrt{-2S-6\delta}$. Since $\alpha\in(0,\sqrt
{-2S} )$,
inequality~\eqref{eqfindoptimalpathtomiddle} will
be violated for large values of $k$ if we choose sufficiently small
$\delta$.

We conclude that with this choice of $\delta$ and $h$,
for sufficiently large $y$,
no path $\gamma\in M_{0,\R,V}^{t,x}$ with $\gamma(0)=y$ can have
$\tau>\tau_{h,\delta}$. On the other hand, if $\tau\le
\tau_{h,\delta}$, then
\[
\Ac_V^{0,\tau}(\gamma)\ge\frac{(k-h)^2}{2\tau_{h,\delta
}}-N\bigl([0,
\tau_{h,\delta}]\times\R\bigr)>V(h)
\]
for sufficiently large $k$, and such a path cannot be a minimizer since
$V(h)$ is the complete action on $[0,\tau]$ for
the
trajectory staying at $h$.

The case of $y\in[-k-1,-k]$ is treated similarly.
\end{pf}

\begin{pf*}{Proof of part~(\hyperlink{itattractor}{5}) of Theorem~\ref
{thmain}} Let
us take any two initial conditions
$v_1=V'_1,v_2=V'_2$ such
that $V_1$ and $V_2$ satisfy~\eqref{eqasymptoticslope}. Then there
is $\alpha\in(0,q)$ such that~\eqref{eqVforlargey} holds
for $V=V_1$ and $V=V_2$.

Let us take $R>L$ given by Lemma~\ref{lmcollapse-with-positive-prob}
and $R_0=R_0(\omega)$ given by
Lemma~\ref{lminitialpointloc}. Due to Lemma~\ref
{lmcollapse-with-positive-prob},
$\Pp\{r^{\pm}<R\}>0$, where $r^\pm$ was introduced in Lemma~\ref
{lmmain-loc}. That lemma, along with the ergodicity of\vadjust{\goodbreak}
the flow $(\theta^t)$
and Poincar\'e Recurrence theorem, allows us to find $n>0$
such that $r^{\pm}(\theta^{-n}\omega)<R$ and \mbox{$\tau^{\pm}=\tau
^{\pm}_{R_0}(\theta^{-n}\omega)<n$}.

If $V=V_1$ or $V=V_2$, then for any $x\in B_{R_0}$ and for sufficiently
large $t$, any $\gamma\in M_{0,\R,V}^{t,x}$
must (by Lemma~\ref{lminitialpointloc}) belong to $M_{0,y}^{t,x}$
for some $y\in B_{R_0}$, and, consequently, Lemma~\ref{lmmain-loc}
implies $\gamma(n)\in B_R$.

Lemma~\ref{lmcollapse-with-positive-prob} and the Poincar\'e
Recurrence theorem imply that there is $n'>n$ and a point
$(t^*,x^*)$
such that for sufficiently large $t$ and for all $z,x\in B_R$, every
$\gamma\in M_{n,z}^{t,x}$ passes through
$(t^*,x^*)$.
Therefore, for these values of $t$ and any $x\in B_L\subset B_R$, any
two minimizers $\gamma_1\in M_{0,\R,V_1}^{t,x}$
and $\gamma_2\in M_{0,\R,V_2}^{t,x}$
pass through $(t^*,x^*)$. Therefore
\[
M_{0,\R,V_1}^{t,x} |_{[t^*,t]}=M_{0,\R,V_2}^{t,x}
|_{[t^*,t]},
\]
which implies
\[
\Phi^{0,t}_{\omega} v_1 |_{B_L}=
\Phi^{0,t}_{\omega} v_2 |_{B_L},
\]
and the forward attraction follows since one can take $v_2=u_\omega$.

The proof of the backward attraction is similar, and we omit it.

The global solution uniqueness also follows automatically if~\eqref
{equniqunessclass} holds with probability 1.
If all we know is that~\eqref{equniqunessclass} holds with positive
probability, then we can use its invariance under
the dynamics and
the ergodicity of $(\theta^t)$ to see that
then it holds with probability 1.
\end{pf*}


\section{Global minimizers}\label{secglobalmin}

A path $\gamma\dvtx\R\to\R$ is called a global minimizer if
$\gamma|_{[t_0,t_1]}\in M_{t_0,\gamma(t_0)}^{t_1,\gamma(t_1)}$ for
any $t_0,t_1$ satisfying $t_0<t_1$.
A global minimizer is called recurrent if there is $R>0$ and a
two-sided sequence $(t_k)_{k\in\Z}$
such that $\lim_{k\to\pm\infty} t_k=\pm\infty$ and $\gamma
(t_k)\in B_R$ for all $k$.

\begin{theorem} With probability 1, there is a unique recurrent global
minimizer $\gamma$.
\end{theorem}

\begin{pf*}{Sketch of proof}
We can use Lemma~\ref{lmcollapse-with-positive-prob} to derive
sequentially that (i) for any $T_1>0$ and for
sufficiently large values $T_2$
the restrictions onto $[-T_1,T_1]$ of minimizers in
$M_{-T_2,0}^{T_2,0}$ can be bounded\vspace*{2pt} by the process of localization
radii $r^{\pm}(\theta^t\omega)$, $t\in[-T_1,T_1]$ ;
(ii) for any $T_1>0$, and for sufficiently large values $T_2$, these
minimizers pass through common points before $-T_1$
and after $T_1$, and, therefore, coincide between these points. We
conclude that as $T_2\to\infty$ the minimizers
stabilize on any finite interval,
and the restriction of the resulting limiting trajectory $\gamma
=\gamma_\omega$ on any finite time interval is a
minimizer. Moreover,
$|\gamma_\omega(t)|<r^{\pm}(\theta^t\omega)$ and the recurrence
property of $\gamma_\omega$ follows.

If $\tilde\gamma$ is another recurrent global minimizer, then again
one can use
Lemma~\ref{lmcollapse-with-positive-prob} to prove that
there is a sequence of times $(s_k)_{k\in\Z}$ such that $\lim_{k\to
\pm\infty} s_k=\pm\infty$ and
$\tilde\gamma(s_k)=\gamma_\omega(s_k)$.
Therefore, $\tilde\gamma$ has to coincide with $\gamma_\omega$.\vadjust{\goodbreak}
\end{pf*}

The following statement can be proven in a similar way:

\begin{theorem} \label{thhyperhyperbolicity}
If $\tilde\gamma$ is one of the one-sided infinite minimizers
constructed in Remark~\ref{remone-sidedminimizers},
then there is $\tau>0$
such that restrictions of $\gamma$ and $\tilde\gamma$ on $(-\infty
,\tau]$ coincide.
\end{theorem}

This property shows that the global minimizer has superstrong
attraction property in the reverse time. In
previously considered situations the
exponential convergence of one-sided minimizers in the reverse time was
a manifestation of hyperbolicity of the global
minimizer.
In analogy with that case, it is natural to refer to the property of
finite-time supercontraction described in
Theorem~\ref{thhyperhyperbolicity}
as ``hyperhyperbolicity.''

\textit{Global minimizers for the Burgers equation with spatially
periodic forcing}.
Let us now change the framework and switch to the Burgers equation with
spatially periodic random forcing.
One of the questions that has not been answered for the periodic
Burgers equation
is the fluctuations of the global minimizer. This question was posed to
the author by Yakov Sinai.
The goal of this section is to prove that for the Poissonian forcing on
the circle $\Sb^1$, the unique global minimizer
has diffusive
behavior.

Let us consider the Burgers dynamics on $\Sb^1=\R/\Z$ under
Poissonian forcing with intensity measure given by $dt\times
m(dx)$ on the
cylinder $\R\times\Sb^1$ for some Borel measure $m$ on $\Sb^1$. If
we restrict ourselves to the set
$\Ub_0=\{u\dvtx\int_{\Sb^1}u=0\}$, then the potential
$V$ of any $v\in\Ub_0$ is well-defined, and we can define the dynamics
via~\eqref{eqdefofcocycle}, \eqref{eqminimizationproblem} with
the only
correction that paths are in $\Sb^1$.

The set $\Ub_0$ is invariant under this dynamics. In fact, between
occurences of Poissonian points, we are solving the
usual
unforced Burgers equation and the mean velocity stays constant. On the
other hand it is easy to see that
the mean velocity is continuous in time even at the time corresponding
to the occurence of a forcing point. A continuous
piecewise constant function is constant, and our invariance claim follows.

The theory that was developed above for the Poisson forcing on the
line, applies to this case as well, so one has
a unique global attracting soliution. One also has a unique global
minimizer $\gamma_\omega$ with asymptotic slope 0
corresponding to the mean
velocity 0. To formulate the main theorem we must unfold $\Sb^1$ onto
its universal cover $\R$ and treat $\gamma_\omega$
as a
continuous path on~$\R$.

\begin{theorem} There is a nonrandom number $D>0$ such that $\gamma
_\omega(t)/\break\sqrt{D|t|}$ converges in distribution
to the standard Gaussian random variable as $t\to\infty$.
\end{theorem}
\begin{pf}
The times between occurrences of Poisson
points are exponentially
distributed. Therefore, the set of all
Poisson points $(t_k(\omega),x_k(\omega))$ such that there\vadjust{\goodbreak}
are no other Poisson points in $[t_k-1,t_k+1]\times\Sb^1$ is
unbounded in both directions.
We agree that $\cdots<t_{-2}<t_{-1}<0<t_{0}<t_1<t_2<\cdots$.
It is easy to check that the global minimizer $\gamma$ passes through
all these points on the cylinder (or their lifts
on
the universal cover).

Let us denote $\Delta_k t=t_k-t_{k-1}$, $\Delta_k x=x_{k}-x_{k-1}
(\mathrm{mod}\ 1)$
and $\Delta_k \gamma= \gamma(t_k)-\gamma(t_{k-1})$.
Notice that all random variables from sequences $(x_k)_{k\in\Z}$,
$(\Delta_k t)_{k\in\Z}$, and realizations of
Poissonian
point field between $t_{k-1}$ and $t_k$ are jointly independent. They
are also identically distributed within each
sequence,
the tails of $\Delta_k t$ are exponential and the distribution of
$x_k$ is $m(dx)/m(\Sb^1)$.

Since
$\Delta_k \gamma$ is a functional of $\Delta_k t$, $\Delta_k x$,
and the realization of the Poissonian field between
$t_{k-1}$
and $t_k$, the sequence $(\Delta_k\gamma)_{k\in\Z}$ of identically
distributed random variables is 1-dependent (the
dependence
comes only through $x_k$ occuring in both $\Delta_k x$ and $\Delta_{k+1}x$).

We know that there is no systematic drift, that is, $(\gamma
(t_k)-\gamma(t_0))/(t_k-t_0)\to0$ as $k\to\infty$.
By the law of large numbers,
$(t_k-t_0)/k\to h =\E(t_1-t_0)$, so $(\gamma(t_k)-\gamma(t_0))/k\to
0$, and $\E(\gamma(t_1)-\gamma(t_0))=0$.

Therefore, by Bernstein's CLT for $m$-dependent random variables, we
conclude that
the distribution of $(\gamma(t_k)-\gamma(t_0))/\sqrt{\sigma^2 k}$
converges weakly to the standard Gaussian one,
where $\sigma^2=\E(\gamma(t_1)-\gamma(t_0))^2 $. Applying the law
of large numbers once again, we conclude that
\[
\frac{\gamma(t)}{\sqrt{{\sigma^2t}/{h} }}\stackrel{d}
{\longrightarrow}\Nc(0,1)
\]
as $t\to\infty$ along the sequence $(t_k)$. To finish the proof, one
has to extend this convergence to all intermediate
values of $t$, but this is not hard since the tails of $\Delta_k t$
are exponential.
This completes the proof with $D=\sigma^2/h$.
\end{pf}

%
\begin{remark} It is also possible to prove a functional version of
the above CLT with two-sided Wiener measure in
the role of
the limiting distribution for appropriately normalized global minimizer.
\end{remark}

\section*{Acknowledgment}
The author thanks Konstantin Khanin for stimulating discussions.

%

%


\printaddresses

\end{document}